\RequirePackage{fix-cm}
\documentclass[smallextended]{svjour3}       
\smartqed  
\usepackage{graphicx}
\usepackage{fullpage}
\usepackage{subfigure}
\begin{document}

\title{Numerical solution of time-fractional Burgers equation in reproducing kernel space
}


\author{Mehmet Giyas Sakar\and Onur Sald{\i}r \and Fevzi Erdogan }


\institute{Mehmet Giyas Sakar \at
              Yuzuncu Yil University, Faculty of Sciences, Department of
Mathematics, 65080, Van, Turkey \\
              Tel.: +90 (432) 2251701\\
              \email{giyassakar@hotmail.com}           
           \and
          Onur Sald{\i}r \at
              Yuzuncu Yil University, Faculty of Sciences, Department of
Mathematics, 65080, Van, Turkey\\
              \email{onursaldir@gmail.com}           
           \and
           Fevzi Erdogan \at
             Yuzuncu Yil University, Faculty of Sciences, Department of
Statistics, 65080, Van, Turkey\\
              \email{ferdogan@yyu.edu.tr}           
}

\date{Received: date / Accepted: date}

\maketitle

\begin{abstract}
In this paper, we present an iterative reproducing kernel method for numerical solution of one dimensional fractional Burgers equation with variable coefficient. Convergence analysis is constructed theoretically. Numerical experiments show that approximate solution uniformly converges to exact solution. The results demonstrate that the given method very efficient and convenient for fractional Burgers equation.
\end{abstract}

\section{Introduction}
Recently, a great deal of important phenomena in dynamical systems, diffusion wave, heat conduction,
cellular systems, oil industries, fluid mechanics, control theory, signal processing and other subjects of sciences
and engineering can be achievely modeled by using differential equations with fractional order (Podlubny 1999). Applications, methods to find approximate solution and qualitative behaviors of solution for fractional
differential equation have been investigated by authors (Diethelm 2010; Lakshmikantham et al. 2009) and references therein.

In this research, an iterative reproducing kernel method is introduced for numerical approach of the time-fractional Burgers equation as follow (Li and Cui 2009):
\begin{eqnarray}
&D_{\eta}^{\alpha }y(\xi,\eta)+k_{1}(\xi,\eta)y_{\xi\xi}(\xi,\eta)+k_{2}(\xi,\eta)y(\xi,\eta)+k_{3}(\xi,\eta)y_{\xi}(\xi,\eta)+k_{4}(\xi,\eta)y(\xi,\eta)y_{\xi}(\xi,\eta)=f(\xi,\eta)\\
&0\leq \xi\leq 1, 0\leq \eta\leq 1, \,\ 0<\alpha\leq 1, \nonumber
\end{eqnarray}
where $k_{1}(\xi,\eta)$, $k_{2}(\xi,\eta)$, $k_{3}(\xi,\eta)$, $k_{4}(\xi,\eta)$ and $f(\xi,\eta)$ are continuous functions. Also $D_{\eta}^{\alpha }$ is Caputo fractional differential operator with respect to time variable $\eta$, subject to
initial and boundary conditions as
\begin{eqnarray}
\left\{
  \begin{array}{ll}
    y(\xi,0)=0, & \\
   y(0,\eta)=y(1,\eta)=0. &
  \end{array}
\right.
\end{eqnarray}
As it well known that many biological, physical and chemical problems are qualified by the transaction of convection and diffusion. Burgers equation is regard as a model problem which determines transaction of convection and diffusion. In recent years, many authors have studied Burgers equation by using several methods. Solution of modified Burgers equation researched with using collocation method by (Ramadan and El-Danaf 2005). A mixed method is presented for one dimensional Burgers equation by (Bahad{\i}r and Sa\u{g}lam 2005). A numerical approximation suggested for Burgers equation with cubic B-splines method by (Da\u{g} et al. 2005). Finite element approximation used for Burgers equation by (Caldwell et al. 1981). Finite difference method is used for numerical approximation of one-dimensional Burgers equation by (Kutluay et al. 1999).

The notion of reproducing kernel goes back to the paper (Zaremba 1908). His work focused on boundary value problems with Dirichlet condition
about harmonic and biharmonic functions. The reproducing property of kernel
functions has a significant role in RKHS theory. Reproducing kernel method
(RKM) gives the solution of differential equations in convergent series form. In the recent times, reproducing
kernel and other methods are carried out for some type of differential and partial differential
equations. For instance, Fredholm integro-differential equations (Arqub et al. 2013), for nonlinear coupled Burgers equations (Mohammadi et al. 2013), eighth order boundary value
problems (Akram and Rehman 2013), fractional Riccati differential equations (Sakar et al. 2017), Riccati differential equations (Sakar 2017),
one-dimensional sine-Gordon equation equations (Akg{\"{u}}l et al. 2016), nonlinear system of PDEs (Mohammadi and Mokhtari 2014), fractional advection-dispersion equation (Jiang and Lin 2010), time
fractional telegraph equation (Jiang and Lin 2011), nonlinear hyperbolic telegraph equation (Yao 2011), reaction-diffusion equations (Lin and Zhou 2004), class of fractional partial
differential equation (Wang et al. 2013) and so on (Sakar and Sald{\i}r 2017; Mohammadi et al. 2016; Cui and Lin 2009).

This study is structured as follows: some definitions and reproducing kernel Hilbert spaces (RKHS) are
presented in Section 2. Representation solution of model problem is illustrated by iterative RKHS method in
Section 3. Convergence of proposed iterative approach is given in Section 4. Some numerical applications of time-fractional Burgers equation are achievely
solved by the introduced method in Section 5. Finally, a conclusion is given in Section 6.
\section{Basic definitions and some reproducing kernel Hilbert space}
\noindent In this part, we provide some important definitions and reproducing kernel Hilbert space which shall be used
in this study.\\\\
\textbf{Definition 2.1} The Caputo  $\alpha$ order derivative (Podlubny 1999; Diethelm 2010) defined as
\begin{eqnarray*}
{\ }^cD_{\eta}^{\alpha}y(\xi,\eta)=\frac{1}{\Gamma(n-\alpha)}\int_{0}^{\eta}\frac{\partial_ry(\xi,r)}{(\eta-r)^{1+\alpha-n}}dr
\end{eqnarray*}
here $n-1<\alpha\leq n$ and $\eta>0$.\\\\
\textbf{Definition 2.2} Let $B\neq\emptyset$  an abstract set. $S:B\times B \to {C}$ function is a reproducing kernel of the Hilbert space
$H$ iff
\begin{eqnarray*}
&&i. \,\,\ \forall \tau \in B, \,\,S(.,\tau)\in H\\
&&ii. \,\,\ \forall \tau \in B, \,\,\forall \mu \in H,\,\langle \mu(.),S(.,\tau)\rangle=\mu(\tau)
\end{eqnarray*}
The second condition (ii) is ``the reproducing property''. Then, this Hilbert space is called RKHS.\\\\

\noindent\textbf{$\mathbf{W_{2}^{1}[0,1]}$ reproducing kernel Hilbert Space}\\\\
$W_{2}^{1} \left[ {0,\,1} \right]$ space is described as
\begin{eqnarray*}
W_{2}^{1}[0,1]=\{g(\xi)|g \,\ \hbox{is absolutely continuous function}, \,\ {g}' \in L^2[0,1]\}.
\end{eqnarray*}
The inner product of $W_{2}^{1} \left[ {0,\,1} \right]$ is
expressed as follow:
\begin{equation}
\langle g\left( \xi \right),f\left( \xi \right)\rangle _{W_{2}^{1}[0,1]} = g\left( 0
\right)f\left( 0 \right) + \int\limits_0^1 {{g}'\left( \xi \right){f}'\left( \xi \right)d\xi},
\end{equation}
and norm of $W_{2}^{1} \left[ {0,\,1} \right]$ is
given as follow:
\[
\left\| g \right\|_{W_{2}^{1}}^2 =\langle g,g\rangle_{W_{2}^{1}}
,\,\,\,g,f \in W_{2}^{1} \left[ {0,\,1} \right].
\]
In [23], the authors showed that $W_{2}^{1} \left[ {0,\,1} \right]$ is
a complete RKHS and its reproducing kernel is,
\begin{equation}
R_x^{\{1\}} \left( \xi \right) = \left\{ {{\begin{array}{*{20}c}
 {1 + \xi , \,\,\, \xi \le x,} \hfill \\
 {1 + x , \,\,\, x >\xi.} \hfill \\
\end{array} }} \right.
\end{equation}
\textbf{$\mathbf{W_{2}^{2}[0,1]}$ reproducing kernel Hilbert Space}\\\\
$W_{2}^{2} \left[ {0,\,1} \right]$ space is defined as;
$W_{2}^{2} \left[ {0,\,1} \right] = \left\{ {\left. {\,g\left( \eta \right)\,}
\right|\,g,\,{g}' \,} \right.$are absolutely continuous functions, $\left. {{g}'' \in L^2\left[
{0,\,1} \right]}, g(0)=0 \right\} $.\\\\
Here, $L^2\left[ {0,1} \right] = \{\left. {g\,} \right|\,\int\limits_0^1
{g^2\left( \eta \right)d\eta < \infty } \}$. The inner product of $W_{2}^{2} \left[ {0,\,1} \right]$  as follow:
\begin{equation}
\langle g\left( \eta \right),f\left( \eta \right)\rangle _{W_{2}^{2}[0,1] } = g\left( 0
\right)f\left( 0 \right) + {g}'\left( 0 \right){f}'\left( 0 \right) +
\int\limits_0^1 {{g}''\left( \eta \right){f}''\left( \eta \right)d\eta}
\end{equation}
and norm of $W_{2}^{2}[0,1]$ given as:
\begin{eqnarray*}
\left\| g \right\|_{W_{2}^{2}}^2 =\langle g,g\rangle_{W_{2}^{2} }
,\,\,\,f,g \in W_{2}^{2} \left[ {0,\,1} \right].
\end{eqnarray*}
\begin{eqnarray}
R_t^{\{2\}} \left( \eta \right) = \left\{ {{\begin{array}{*{20}c}
 {\eta t+\frac{1}{2}t \eta^{2}-\frac{1}{6}\eta^{3} , \,\,\,\eta \le t,} \hfill \\
 {-\frac{1}{6}t^{3}+\frac{1}{2}\eta t^{2}+t\eta , \,\,\,\eta > t.} \hfill \\
\end{array} }} \right.
\end{eqnarray}
\textbf{$\mathbf{W_{2}^{3}[0,1]}$ reproducing kernel Hilbert Space}\\\\
$W_{2}^{3} \left[ {0,\,1} \right]$ space is defined as;

$W_{2}^{3} \left[ {0,\,1} \right] = \left\{ {\left. {\,g\left( \xi \right)\,}
\right|\,g,\,{g}',\, {g}'' \,} \right.$are absolutely continuous functions, $\left. {{g}^{(3)} \in L^2\left[
{0,\,1} \right]}, g(0)=g(1)=0 \right\} $.\\\\
The inner product of $W_{2}^{3}[0,1]$ as follow:
\begin{equation}
\langle g\left( \xi \right),f\left( \xi \right)\rangle _{W_{2}^{3}[0,1] } = g\left( 0
\right)f\left( 0 \right) + {g}'\left( 0 \right){f}'\left( 0 \right) + g\left( 1 \right)f\left( 1 \right) +
\int\limits_0^1 {{g}^{(3)}\left( \xi \right){f}^{(3)}\left( \xi \right)d\xi}
\end{equation}
and norm of $W_{2}^{3}[0,1]$ given as:
\[
\left\| g \right\|_{W_{2}^{3}}^2 = \langle g,g\rangle_{W_{2}^{3} }
,\,\,\,f,g \in W_{2}^{3} \left[ {0,\,1} \right].
\]

\begin{equation}
R_x^{\{3\}}(\xi) = \left\{ {{\begin{array}{*{20}c}
 {\frac{-1}{120}(x-1)\xi(\xi x^4-4\xi x^3+6\xi x^2+x\xi^4-5x\xi^3-120x\xi+120x+\xi^4),\,\,\,\xi \le x,} \hfill \\\\
 {\frac{-1}{120}(\xi-1)x(x\xi^4-4x\xi^3+6x\xi^2+\xi x^4-5\xi x^3-120\xi x +120\xi+x^4),\,\,\,\xi > x,} \hfill \\
\end{array} }} \right.
\end{equation}

\noindent\textbf{$\mathbf{W_{2}^{(3,2)}(D)}$ reproducing kernel Hilbert Space}\\\\
Let be $D=[0,1]\times[0,1]$. $W_{2}^{(3,2)}(D)$ space defined as follow;
\begin{eqnarray}
W_{2}^{(3,2)}(D)=\{y(\xi,\eta)|\frac{\partial^3y}{\partial \xi^2\partial \eta}\,\ \hbox{is completely continuous in}\,\ D,\nonumber\\
 \frac{\partial^5y}{\partial \xi^3\partial \eta^2}\in L^2(D),y(\xi,0)=y(0,\eta)=y(1,\eta)=0\}
\end{eqnarray}
and the inner product given as:
\begin{eqnarray}
&&\langle y(\xi,\eta),v(\xi,\eta) \rangle_{W_{2}^{(3,2)}}=\sum\limits_{i=0}^{1}\int\limits_0^1[\frac{\partial^2}{\partial \eta^2}\frac{\partial^i}{\partial \xi^i}y(0,\eta)\frac{\partial^2}{\partial \eta^2}\frac{\partial^i}{\partial \xi^i}v(0,\eta)]d\eta+\int_{0}^{1}\frac{\partial^2}{\partial \eta^2}y(1,\eta)\frac{\partial^2}{\partial \eta^2}v(1,\eta)]d\eta\nonumber\\
&+&\sum_{j=0}^{1}\langle \frac{\partial^j}{\partial \eta^j}y(\xi,0),\frac{\partial^j}{\partial \eta^j}v(\xi,0) \rangle_{W_2^3}
+\int\limits_0^1\int\limits_0^1[\frac{\partial^3}{\partial \xi^3}\frac{\partial^2}{\partial \eta^2}y(\xi,\eta)\frac{\partial^3}{\partial \xi^3}\frac{\partial^2}{\partial \eta^2}v(\xi,\eta)]d\xi d\eta
\end{eqnarray}
and norm of $W_{2}^{(3,2)}(D)$ given as:
\begin{eqnarray*}
\left\| y \right\|_{W_{2}^{(3,2)}}^2 = \langle y,y\rangle_{W_{2}^{(3,2)}}
,\,\,\,y,v \in W_{2}^{(3,2)}(D).
\end{eqnarray*}
\textbf{Theorem 2.2} Let $K_{(x,t)}(\xi,\eta)$ be a reproducing kernel of $W_{2}^{(3,2)}(D)$. So, we can write
\begin{eqnarray*}
K_{(x,t)}(\xi,\eta)=R_x^{\{3\}}(\xi)R_t^{\{2\}}(\eta).
\end{eqnarray*}
where $R_x^{\{3\}}(\xi)$ and $R_t^{\{2\}}(\eta)$ are reproducing kernel functions of $W_2^3[0,1]$ and $W_2^2[0,1]$, respectively. For any $y(\xi,\eta)\in W_{2}^{(3,2)}(D)$
\begin{eqnarray*}
y(x,t)=\langle y(\xi,\eta), K_{(x,t)}(\xi,\eta)\rangle_{W_2^{(3,2)}}
\end{eqnarray*}
and
\begin{eqnarray}
K_{(\xi,\eta)}(x,t)=K_{(x,t)}(\xi,\eta)
\end{eqnarray}
\textbf{$\mathbf{W_{2}^{(1,1)}(D)}$ reproducing kernel Hilbert Space}\\\\
$W_{2}^{(1,1)}(D)$ space is given as,
\begin{eqnarray*}
W_{2}^{(1,1)}(D)=\{y(\xi,\eta)| \,\ y \,\ \hbox{is completely continuous in}\,\ D=[0,1]\times[0,1], \frac{\partial^2y}{\partial \xi\partial \eta}\in L^2(D) \}
\end{eqnarray*}
The inner product of $W_{2}^{(1,1)}(D)$ as follow:
\begin{eqnarray*}
\langle y(\xi,\eta),v(\xi,\eta) \rangle_{W_{2}^{(1,1)}}&=&\int\limits_0^1[\frac{\partial}{\partial \eta}y(0,\eta)\frac{\partial}{\partial \eta}v(0,\eta)]d\eta+\langle y(\xi,0),v(\xi,0) \rangle_{W_2^1}\\
&+&\int\limits_0^1\int\limits_0^1[\frac{\partial}{\partial \xi}\frac{\partial}{\partial \eta}y(\xi,\eta)\frac{\partial}{\partial \xi}\frac{\partial}{\partial \eta}v(\xi,\eta)]d\xi d\eta
\end{eqnarray*}
and norm defined as:
\begin{eqnarray*}
\left\| y \right\|_{W_{2}^{(1,1)}}^2 =\langle y,y\rangle_{W_{2}^{(1,1)}}
,\,\,\,y,v \in W_{2}^{(1,1)}(D).
\end{eqnarray*}
$W_{2}^{(1,1)}(D)$ space is a RKHS, and its reproducing kernel function $\tilde{K}_{(x,t)}(\xi,\eta)$ is given as
\begin{eqnarray*}
\tilde{K}_{(x,t)}(\xi,\eta)=R_x^{\{1\}}(\xi)R_t^{\{1\}}(\eta)
\end{eqnarray*}
\section {Representation solution of Eqs. (1)-(2) in $W_{2}^{(3,2)}(D)$}
\noindent The representation solution of (1)-(2) will be consisted in $W_{2}^{(3,2)}(D)$. Firstly, we will describe the
linear operator $L$ as,
\begin{equation}
L:W_{2}^{(3,2)}(D) \to W_{2}^{(1,1)}(D),
\end{equation}
such that
\begin{equation}
Ly(\xi,\eta)=D_{\eta}^{\alpha }y+k_{1}(\xi,\eta)y_{\xi\xi}+k_{2}(\xi,\eta)y+k_{3}(\xi,\eta)y_{\xi}
\end{equation}
The problem (1)-(2) can be written as follow:
\begin{equation}
\left\{ {{\begin{array}{*{20}c}
 {Ly\left( \xi,\eta \right) = F(\xi,\eta,y(\xi,\eta),y_{\xi}(\xi,\eta)),\,\,\,\,\,\xi \in \left[
{0,\,1} \right], \,\ \eta\in [0,1]} \hfill \\
 y(\xi,0)=y(0,\eta)=y(1,\eta)=0 \hfill \\
\end{array} }} \right.
\end{equation}
here $F(\xi,\eta,y(\xi,\eta),y_{\xi}(\xi,\eta))= f(\xi,\eta)-k_{4}(\xi,\eta)y(\xi,\eta)y_{\xi}(\xi,\eta)$.\\\\
Let we choose a countable dense subset $\{(\xi_i,\eta_i)\}_{i=1}^{\infty} \in D$, defining
\begin{eqnarray}\label{eq9}
\psi_i(\xi,\eta)&=&L_{(x,t)}K_{(x,t)}(\xi,\eta)|_{(x,t)=(\xi_{i},\eta_{i})}\nonumber\\
&=&\{D_{t}^{\alpha}K_{(x,t)}(\xi,\eta)+k_{1}(x,t)\frac{\partial^2}{\partial x^2}K_{(x,t)}(\xi,\eta)+k_{2}(x,t)K_{(x,t)}(\xi,\eta)\nonumber\\
&+&k_{3}(x,t)\frac{\partial}{\partial x}K_{(x,t)}(\xi,\eta)\}|_{(x,t)=(\xi_{i},\eta_{i})}\nonumber\\
&=&\frac{1}{\Gamma(1-\alpha)}\int_{0}^{\eta_{i}}\frac{\partial_{r}K_{(\xi,r)}(\xi,\eta)}{(\eta_{i}-r)^{\alpha}}dr+k_{1}(\xi_i,\eta_i)\frac{\partial^2}{\partial x^2}K_{(\xi_i,\eta_i)}(\xi,\eta)+k_{2}(\xi_i,\eta_i)K_{(\xi_i,\eta_i)}(\xi,\eta)\nonumber\\
&+&k_{3}(\xi_i,\eta_i)\frac{\partial}{\partial x}K_{(\xi_i,\eta_i)}(\xi,\eta), i=1,2,...
\end{eqnarray}
where $K_{(x,t)}(\xi,\eta)$ is the reproducing kernel of $W_{2}^{(3,2)}(D)$.\\\\
\textbf{Theorem 3.1} $\psi_i(\xi,\eta)\in W_{2}^{(3,2)}(D), i=1,2,...$\\\\
\textbf{Proof.} Using definition of $W_{2}^{(3,2)}(D)$ space, firstly we will show that $\frac{\partial^5\psi_i(\xi,\eta)}{\partial\xi^3\partial\eta^2}\in L^2(D)$ and $\frac{\partial^3\psi_i(\xi,\eta)}{\partial\xi^2\partial\eta}$ is completely continuous function. Then, we will demonstrate that $\psi_i(\xi,\eta)$ satisfies the initial and boundary conditions.\\
Now, by property of kernel function $K_{(x,t)}(\xi,\eta)$ we can take,
\begin{eqnarray*}
\partial_{x^2\xi^3\eta^2}^{7}K_{(x,t)}(\xi,\eta)=\partial_{x^2\xi^3}^{5}R_{x}^{\{3\}}(\xi)\partial_{\eta^2}^{2}R_{t}^{\{2\}}(\eta)
\end{eqnarray*}
In here, $\partial_{x^2\xi^3}^{5}R_{x}^{\{3\}}(\xi)$ and $\partial_{\eta^2}^{2}R_{t}^{\{2\}}(\eta)$ are continuous functions on $[0,1]$. Because these functions are continuous on closed interval, so these functions are bounded. Therefore we can write,
\begin{eqnarray*}
|\partial_{x^2\xi^3\eta^2}^{7}K_{(x,t)}(\xi,\eta)|\leq M_1
\end{eqnarray*}
In a similar manner, one can see that
\begin{eqnarray*}
|\partial_{t\xi^3\eta^2}^{6}K_{(x,t)}(\xi,\eta)|\leq M_2\\
|\partial_{\xi^3\eta^2}^{5}K_{(x,t)}(\xi,\eta)|\leq M_3\\
|\partial_{x\xi^3\eta^2}^{6}K_{(x,t)}(\xi,\eta)|\leq M_4
\end{eqnarray*}
Here, $M_1,M_2,M_3$ and $M_4$ are positive constants. From (\ref{eq9}),
\begin{eqnarray*}
|\frac{\partial^5\psi_i(\xi,\eta)}{\partial\xi^3\partial\eta^2}|&\leq&|\frac{1}{\Gamma(1-\alpha)}\int_{0}^{\eta_i}\frac{M_2}{(\eta_i-r)^{\alpha}}dr+k_{1}(\xi_i,\eta_i)M_1\\
&+&k_{2}(\xi_i,\eta_i)M_3+k_{3}(\xi_i,\eta_i)M_4|\\
&\leq&\frac{M_2}{\Gamma(2-\alpha)}\eta_i^{1-\alpha}+|k_1(\xi_i,\eta_i)|M_1+|k_2(\xi_i,\eta_i)|M_3+|k_3(\xi_i,\eta_i)|M_4
\end{eqnarray*}
Therefore, $\frac{\partial^5\psi_i(\xi,\eta)}{\partial\xi^3\partial\eta^2}\in L^2(D)$. Noting that $D$ is closed, thus, $\frac{\partial^3\psi_i(\xi,\eta)}{\partial\xi^2\partial\eta}$ is completely continuous in $D$. $\psi_i(\xi,\eta)$ satisfies the conditions of problem so that $K_{(x,t)}(\xi,0)=0$ and $K_{(x,t)}(0,\eta)=K_{(x,t)}(1,\eta)=0$. Thus $\psi_i(\xi,\eta)\in W_{2}^{(3,2)}(D)$. The proof is completed.\\\\
\textbf{Theorem 3.2} $\{\psi_i(\xi,\eta)\}_{i=1}^{\infty}$ is a complete system in $W_{2}^{(3,2)}(D)$.\\\\
\textbf{Proof.} We have
\begin{eqnarray}
\psi_i(\xi,\eta)&=&(L^{\ast}\varphi_i)(\xi,\eta)=\langle (L^{\ast}\varphi_i)(x,t),K_{(\xi,\eta)}(x,t) \rangle_{W_{2}^{(3,2)}}\nonumber\\
&=&\langle \varphi_i(x,t),L_{(x,t)}K_{(\xi,\eta)}(x,t) \rangle_{W_{2}^{(1,1)}}=L_{(x,t)}K_{(\xi,\eta)}(x,t)|_{(x,t)=(\xi_i,\eta_i)}\nonumber\\
&=&L_{(x,t)}K_{(x,t)}(\xi,\eta)|_{(x,t)=(\xi_i,\eta_i)}.
\end{eqnarray}
Clearly, $\psi_i(\xi,\eta)\in W_{2}^{(3,2)}(D)$, for each fixed $y(\xi,\eta)\in W_{2}^{(3,2)}(D)$, if $\langle y(\xi,\eta),\psi_i(\xi,\eta) \rangle_{W_{2}^{(3,2)}}=0, \,\,\ i=1,2,...$.
So,
\begin{eqnarray}
\langle y(\xi,\eta),(L^{\ast}\varphi_i)(\xi,\eta) \rangle_{W_{2}^{(3,2)}}=\langle Ly(\xi,\eta),\varphi_i(\xi,\eta) \rangle_{W_{2}^{(1,1)}}=(Ly)(\xi_i,\eta_i)=0, \,\,\ i=1,2,...
\end{eqnarray}
$(Ly)(\xi,\eta)=0$ since $\{(\xi_i,\eta_i)\}_{i=1}^{\infty}$ is dense in $D$. Thus, $y=0$ by the existence of $L^{-1}$. The proof is completed.\\\\

The orthonormal system $\{\overline{\psi_i}(\xi,\eta)\}_{i=1}^{\infty}$ of $W_{2}^{(3,2)}(D)$ can be obtained by the Gram-Schmidt orthogonalization of $\{\psi_i(\xi,\eta)\}_{i=1}^{\infty}$ as
\begin{eqnarray}
\overline{\psi_i}(\xi,\eta)=\sum\limits_{k=1}^i \beta_{ik}\psi_{k}(\xi,\eta).
\end{eqnarray}
Here, $\beta_{ik}$ are orthogonalization coefficients, $\beta_{ii}>0$, $i=1,2,...$\\\\
\textbf{Theorem 3.3} If $\{(\xi_i,\eta_i)\}_{i=1}^{\infty}$ is dense in $D$, then the solution (14) is
\begin{eqnarray}
y(\xi,\eta)=\sum\limits_{i=1}^{\infty}\sum\limits_{k=1}^{i}\beta_{ik}F(\xi_k,\eta_k,y(\xi_k,\eta_k),\partial_{\xi}y(\xi_k,\eta_k))\overline{\psi}_{i}(\xi,\eta)
\end{eqnarray}
\textbf{Proof.}$\{\psi_i(\xi,\eta)\}_{i=1}^{\infty}$ is complete system in $W_{2}^{(3,2)}(D)$. Therefore, we get
\begin{eqnarray}
y(\xi,\eta)&=&\sum\limits_{i=1}^{\infty}\langle y(\xi,\eta),\overline{\psi}_i(\xi,\eta \rangle_{W_{2}^{(3,2)}}\overline{\psi}_i(\xi,\eta)=\sum\limits_{i=1}^{\infty}\sum\limits_{k=1}^{i}\beta_{ik}\langle y(\xi,\eta),\psi_k(\xi,\eta) \rangle_{W_{2}^{(3,2)}}\overline{\psi}_i(\xi,\eta)\nonumber\\
&=&\sum\limits_{i=1}^{\infty}\sum\limits_{k=1}^{i}\beta_{ik}\langle y(\xi,\eta),L^{\ast}\varphi_k(\xi,\eta) \rangle_{W_{2}^{(3,2)}}\overline{\psi}_i(\xi,\eta)=\sum\limits_{i=1}^{\infty}\sum\limits_{k=1}^{i}\beta_{ik}\langle Ly(\xi,\eta),\varphi_k(\xi,\eta) \rangle_{W_{2}^{(1,1)}}\overline{\psi}_i(\xi,\eta)\nonumber\\
&=&\sum\limits_{i=1}^{\infty}\sum\limits_{k=1}^{i}\beta_{ik}\langle Ly(\xi,\eta),\tilde{K}_{(\xi_k,\eta_k)}(\xi,\eta) \rangle_{W_{2}^{(1,1)}}\overline{\psi}_i(\xi,\eta)=\sum\limits_{i=1}^{\infty}\sum\limits_{k=1}^{i}\beta_{ik}Ly(\xi_k,\eta_k)\overline{\psi}_i(\xi,\eta)\nonumber\\
&=&\sum\limits_{i=1}^{\infty}\sum\limits_{k=1}^{i}\beta_{ik}F(\xi_k,\eta_k,y(\xi_k,\eta_k),\partial_\xi y(\xi_k,\eta_k))\overline{\psi}_i(\xi,\eta)
\end{eqnarray}
So, the proof of theorem is completed.\\
Now an approximate solution $y_n(\xi,\eta)$ can be written by taking finitely $n$-terms of the exact solution $y(\xi,\eta)$ as follows,
\begin{eqnarray}
y_n(\xi,\eta)=\sum\limits_{i=1}^{n}\sum\limits_{k=1}^{i}\beta_{ik}F(\xi_k,\eta_k,y(\xi_k,\eta_k),\partial_\xi y(\xi_k,\eta_k))\overline{\psi}_i(\xi,\eta).
\end{eqnarray}
Its clear that,
\begin{eqnarray}
\|y(\xi,\eta)-y_n(\xi,\eta)\|\rightarrow0 \,\ \hbox{as}\,\ n\rightarrow\infty
\end{eqnarray}
\section{Convergence analysis of iterative process}
If we write
\begin{eqnarray}
A_{i}=\sum\limits_{k=1}^{i}\beta_{ik}F(\xi_k,\eta_k,y(\xi_k,\eta_k),\partial_\xi y(\xi_k,\eta_k)),
\end{eqnarray}
then (19) can be described as
\begin{eqnarray}
y(\xi,\eta)=\sum\limits_{i=1}^{\infty}A_{i}\overline{\psi}_i(\xi,\eta).
\end{eqnarray}
Now we take $(\xi_1,\eta_1)=0$; then from the initial conditions of (14), $y(\xi_1,\eta_1)$ is known. We put $y_0(\xi_1,\eta_1)=y(\xi_1,\eta_1)$ and define the $n$-term approximation to $y(\xi,\eta)$ by
\begin{eqnarray}
y_n(\xi,\eta)=\sum\limits_{i=1}^{n}B_{i}\overline{\psi}_i(\xi,\eta),
\end{eqnarray}
here
\begin{eqnarray}
B_i=\sum\limits_{k=1}^{i}\beta_{ik}F(\xi_k,\eta_k,y_{k-1}(\xi_k,\eta_k),\partial_\xi y_{k-1}(\xi_k,\eta_k)).
\end{eqnarray}
After then, we will confirm that $y_n(\xi,\eta)$ uniformly converges to $y(\xi,\eta)$. Therefore, the following lemma will be given.\\

\noindent\textbf{Lemma 4.1} If $F(\xi,\eta,y(\xi,\eta),y_{\xi}(\xi,\eta))$ is continuous and $y_n\rightarrow \hat{y}$ for $(\xi_n,\eta_n)\rightarrow(x,t)$, then
\begin{eqnarray}
F(\xi_n,\eta_n,y_{n-1}(\xi_n,\eta_n),\partial_{\xi}y_{n-1}(\xi_n,\eta_n))\rightarrow F(x,t,\hat{y}(x,t),\partial_{\xi}\hat{y}(x,t)).
\end{eqnarray}
\textbf{Proof.} Since
\begin{eqnarray}
|y_{n-1}(\xi_n,\eta_n)-\hat{y}(x,t)|&=|y_{n-1}(\xi_n,\eta_n)-y_{n-1}(x,t)+y_{n-1}(x,t)-\hat{y}(x,t)|\nonumber\\
&\leq|y_{n-1}(\xi_n,\eta_n)-y_{n-1}(x,t)|+|y_{n-1}(x,t)-\hat{y}(x,t)|
\end{eqnarray}
From the reproducing kernel feature, we have
\begin{eqnarray}
y_{n-1}(\xi_n,\eta_n)= \langle y_{n-1}(\xi,\eta),K_{(\xi_n,\eta_n)}(\xi,\eta) \rangle_{W_{2}^{(3,2)}}, \,\,\ y_{n-1}(x,t)= \langle y_{n-1}(\xi,\eta),K_{(x,t)}(\xi,\eta) \rangle_{W_{2}^{(3,2)}} .
\end{eqnarray}
It follows that
\begin{eqnarray}
|y_{n-1}(\xi_n,\eta_n)-y_{n-1}(x,t)|=|\langle y_{n-1}(\xi,\eta), K_{(\xi_n,\eta_n)}(\xi,\eta)-K_{(x,t)}(\xi,\eta)  \rangle|.
\end{eqnarray}
From the convergence of $y_{n-1}(\xi,\eta)$, there exists a constant $M$, such that
\begin{eqnarray}
\|y_{n-1}(\xi,\eta)\|_{W_{2}^{(3,2)}} \leq M \|\hat{y}(x,t)\|_{W_{2}^{(3,2)}}, \,\ \hbox{as} \,\ n\geq M.
\end{eqnarray}
As the same time, we can prove
\begin{eqnarray}
\|K_{(\xi_n,\eta_n)}(\xi,\eta)-K_{(x,t)}(\xi,\eta)\|_{W_{2}^{(3,2)}}\rightarrow 0, \,\ \hbox{for} \,\ n\rightarrow \infty
\end{eqnarray}
by using Theorem 2.2. So,
\begin{eqnarray}
y_{n-1}(\xi_n,\eta_n)\rightarrow \hat{y}(x,t), \,\ \hbox{as} \,\ (\xi_n,\eta_n)\rightarrow (x,t).
\end{eqnarray}
In a similar way it can be shown that
\begin{eqnarray}
\partial_{\xi}y_{n-1}(\xi_n,\eta_n)\rightarrow \partial_{\xi}\hat{y}(x,t), \,\ \hbox{as} \,\ (\xi_n,\eta_n)\rightarrow (x,t).
\end{eqnarray}
Therefore,
\begin{eqnarray}
F(\xi_n,\eta_n,y_{n-1}(\xi_n,\eta_n),\partial_{\xi}y_{n-1}(\xi_n,\eta_n))\rightarrow F(x,t,\hat{y}(x,t),\partial_{\xi}\hat{y}(x,t)).
\end{eqnarray}
So, the proof of theorem is completed.\\

\noindent\textbf{Theorem 4.1} Suppose that $\|y_n\|$ is a bounded in (14) and (25) has a unique solution. If $\{(\xi_i,\eta_i)\}_{i=1}^{\infty}$ is dense in $D$, so the $n$- term approximate solution $y_n(\xi,\eta)$ converges to the exact solution $y(\xi,\eta)$ of (14) and
\begin{eqnarray}
y(\xi,\eta)=\sum\limits_{i=1}^{\infty}B_i\overline{\psi}_i(\xi,\eta)
\end{eqnarray}
where $B_i$ is given (26).\\

\noindent\textbf{Proof.} Firstly, we will demonstrate the convergence of $y_n(\xi,\eta)$. From (25), we deduce that
\begin{eqnarray}
y_{n+1}(\xi,\eta)=y_n(\xi,\eta)+B_{n+1}\overline{\psi}_{n+1}(\xi,\eta).
\end{eqnarray}
Based on the orthonormality of $\{\overline{\psi}_{i}\}_{i=1}^{\infty}$ provides that
\begin{eqnarray}
\|y_{n+1}\|^{2}=\|y_n\|^{2}+B_{n+1}^2=\sum\limits_{i=1}^{n+1}B_i^2.
\end{eqnarray}
Therefore, $\|y_{n+1}\|>\|y_n\|$ inequality holds from (38). By reason of the boundedness of $\|y_n\|$, one can see that $\|y_n\|$ is convergent and there exists a constant $c$ such that
\begin{eqnarray}
\sum\limits_{i=1}^{\infty}B_i^2=c.
\end{eqnarray}
So, this indicates that $\{B_i\}_{i=1}^{\infty}\in l^2$. If $m>n$, then
\begin{eqnarray}
\|y_m-y_n\|^2&=&\|y_{m}-y_{m-1}+y_{m-1}-y_{m-2}+\cdots+y_{n+1}-y_{n}\|^2\nonumber\\
&=&\|y_{m}-y_{m-1}\|^2+\|y_{m-1}-y_{m-2}\|^2+\cdots+\|y_{n+1}-y_{n}\|^2.
\end{eqnarray}
On account of
\begin{eqnarray}
\|y_{m}-y_{m-1}\|^2=B_m^2,
\end{eqnarray}
consequently
\begin{eqnarray}
\|y_{m}-y_{n}\|^2=\sum\limits_{l=n+1}^{m}B_l^2\rightarrow 0, \,\ \hbox{as} \,\ n\rightarrow\infty.
\end{eqnarray}
The completeness of $W_{2}^{(3,2)}(D)$ shows that $y_n\rightarrow \hat{y}$ as $n\rightarrow\infty$. Now, we will show that $\hat{y}$ is the representation solution of (14). Taking limits in (25) we get
\begin{eqnarray}
\hat{y}(\xi,\eta)=\sum\limits_{i=1}^{\infty}B_i\overline{\psi}_i(\xi,\eta).
\end{eqnarray}
Note that
\begin{eqnarray}
(L\hat{y})(\xi,\eta)=\sum\limits_{i=1}^{\infty}B_iL\overline{\psi}_i(\xi,\eta),
\end{eqnarray}
\begin{eqnarray}
(L\hat{y})(\xi_l,\eta_l)&=&\sum\limits_{i=1}^{\infty}B_iL\overline{\psi}_i(\xi_l,\eta_l)=\sum\limits_{i=1}^{\infty}B_i\langle L\overline{\psi}_i(\xi,\eta),\varphi_l(\xi,\eta)\rangle_{W_{2}^{(1,1)}}\nonumber\\
&=&\sum\limits_{i=1}^{\infty}B_i\langle \overline{\psi}_i(\xi,\eta),L^{\ast}\varphi_l(\xi,\eta)\rangle_{W_{2}^{(3,2)}}=\sum\limits_{i=1}^{\infty}B_i\langle \overline{\psi}_i(\xi,\eta),\psi_l(\xi,\eta)\rangle_{W_{2}^{(3,2)}}.
\end{eqnarray}
Therefore,
\begin{eqnarray}
\sum\limits_{l=1}^{i}\beta_{il}(L\hat{y})(\xi_l,\eta_l)&=&\sum\limits_{i=1}^{\infty}B_i\langle \overline{\psi}_i(\xi,\eta),\sum\limits_{l=1}^{i}\beta_{il}\psi_l(\xi,\eta)\rangle_{W_{2}^{(3,2)}}\\
&=&\sum\limits_{i=1}^{\infty}B_i\langle \overline{\psi}_i(\xi,\eta),\overline{\psi}_l(\xi,\eta)\rangle_{W_{2}^{(3,2)}}=B_l.
\end{eqnarray}
In view of (38), we have
\begin{eqnarray}
L\hat{y}(\xi_l,\eta_l)=F(\xi_l,\eta_l,y_{l-1}(\xi_l,\eta_l),\partial_{\xi}y_{l-1}(\xi_l,\eta_l)).
\end{eqnarray}
Since $\{(\xi_i,\eta_i)\}_{i=1}^{\infty}$ is dense in $D$, for each $(x,t)\in D$, there exists a subsequence $\{(\xi_{n_j},\eta_{n_j})\}_{j=1}^{\infty}$ such that $(\xi_{n_j},\eta_{n_j})\rightarrow (x,t)$, $j\rightarrow \infty$. We know that
\begin{eqnarray}
L\hat{y}(\xi_{n_j},\eta_{n_j})=F(\xi_{n_j},\eta_{n_j},y_{n_{j-1}}(\xi_{n_j},\eta_{n_j}),\partial_{\xi}y_{n_{j-1}}(\xi_{n_j},\eta_{n_j})).
\end{eqnarray}
Let $j\rightarrow\infty$, the continuity of $F$ and by Lemma 4.1, we can write
\begin{eqnarray}
(L\hat{y})(x,t)=F(x,t,\hat{y}(x,t),\partial_{\xi}\hat{y}(x,t))
\end{eqnarray}
which demonstrates that $\hat{y}(\xi,\eta)$ provide (14). The proof is completed.

\section{Numerical Applications}

In this section, the iterative reproducing kernel approach is implemented for two time-fractional Burgers problems which have exact solution. \\\\
$\mathbf{Example ~5.1}$ We take into consideration the following time-fractional Burgers equation (Cui and Lin 2009):
\begin{equation}
D_{\eta}^{\alpha }y+(1+\xi\eta)y_{\xi\xi}+\xi^2y+(\xi+1)y_{\xi}-\eta\sin(\xi)yy_{\xi}=f(\xi,\eta)
\end{equation}
\begin{eqnarray*}
0\leq \xi\leq 1, 0\leq \eta\leq 1, \,\ 0<\alpha\leq 1,
\end{eqnarray*}
\begin{eqnarray*}
f(\xi,\eta)&=&\frac{(\xi^2-\xi)\eta\pi}{\sin(\pi\alpha)\Gamma(-1-\alpha)}+2(\eta\xi+1)\eta^{1+\alpha}+(\xi^4-\xi^3)\eta^{1+\alpha}\\
&+&(1+\xi)(2\xi-1)\eta^{1+\alpha}-\eta\sin(\xi)(\xi^2-\xi)(\eta^{2+2\alpha})(2\xi-1)
\end{eqnarray*}
initial and boundary conditions of problem as follow:
\begin{eqnarray}
\left\{
  \begin{array}{ll}
    y(\xi,0)=0, & \\
   y(0,\eta)=y(1,\eta)=0, &
  \end{array}
\right.
\end{eqnarray}
The exact solution of problem is given as:
\begin{equation}
y(\xi,\eta) = (\xi^2-\xi)\eta^{1+\alpha}.
\end{equation}
Taking
$\xi_i = \frac{i }{p},\,i = 1,\,2,...,p$, $\eta_i = \frac{i }{q},\,i = 1,\,2,...,q$ and $n=p\times q$. In order to demonstrate of the influence and applicability of the method, the absolute error of Example 5.1 is computed and given in Table 1, Table 2 and Table 3 with different values of $\alpha$ for $n=25 (p=q=5)$. The graphics  of approximate solution $y_n \left( \xi,\eta \right)$ are given in Figure 1 with different values of $\alpha=0.7$, $\alpha=0.8$, and $\alpha=0.9$ respectively.\\

\begin{figure}[ht]
\centerline{\includegraphics[width=2.5000in,height=1.75in]{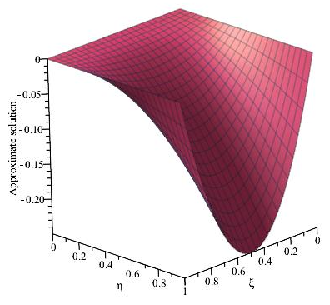}
\includegraphics[width=2.5000in,height=1.75in]{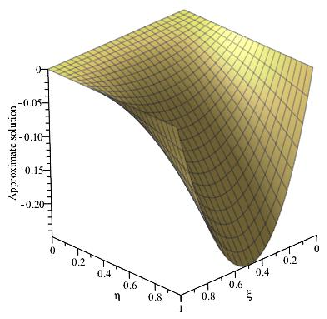}
\includegraphics[width=2.5000in,height=1.75in]{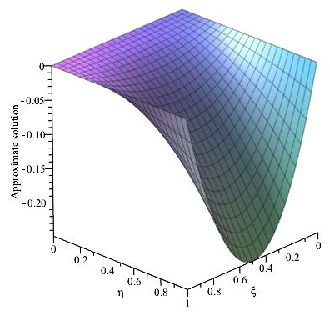}}
\caption{The surface shows approximate solution of Example 5.1 for $n=25$ with $\alpha=0.7$, $\alpha=0.8$ and $\alpha=0.9$ respectively on $0\leq \xi,\eta\leq 1$ }
\label{fig1}
\end{figure}

\begin{table}[h]
\caption{The absolute error of Example 5.1 for $\alpha=0.9$ ($0\leq \xi,\eta\leq 0.6$)}
\begin{tabular}{lllllll}
\hline\noalign{\smallskip}
$\xi/\eta$ & $0.1$ & $0.2$ & $0.3$ & $0.4$ & $0.5$ & $0.6$  \\
\noalign{\smallskip}\hline\noalign{\smallskip}
0.1& 2.39E-4&2.26E-5&3.74E-5&5.31E-5&4.92E-5&7.43E-5 \\
0.2& 4.25E-4&4.23E-5&7.26E-5&1.05E-4&1.05E-4&1.60E-4 \\
0.3& 5.58E-4&5.78E-5&1.01E-4&1.48E-4&1.53E-4&2.33E-4 \\
0.4& 6.37E-4&6.53E-5&1.16E-4&1.71E-4&1.80E-4&2.77E-4 \\
0.5& 6.62E-4&6.74E-5&1.21E-4&1.79E-4&1.89E-4&2.95E-4 \\
0.6& 6.33E-4&6.10E-5&1.11E-4&1.65E-4&1.73E-4&2.77E-4 \\
\noalign{\smallskip}\hline
\end{tabular}
\end{table}

\begin{table}[h]
\caption{The absolute error of Example 5.1 for $\alpha=0.8$ ($0\leq \xi,\eta\leq 0.6$)}
\begin{tabular}{lllllll}
\hline\noalign{\smallskip}
$\xi/\eta$ & $0.1$ & $0.2$ & $0.3$ & $0.4$ & $0.5$ & $0.6$ \\
\noalign{\smallskip}\hline\noalign{\smallskip}
0.1& 2.46E-4&1.09E-5&1.47E-5&1.97E-5&1.14E-5&5.62E-5 \\
0.2& 4.39E-4&1.54E-5&3.50E-5&4.96E-5&4.19E-5&1.31E-4 \\
0.3& 5.78E-4&1.54E-5&5.50E-5&7.80E-5&7.29E-5&1.97E-4 \\
0.4& 6.62E-4&1.46E-5&6.81E-5&9.58E-5&9.22E-5&2.38E-4 \\
0.5& 6.91E-4&1.07E-5&7.71E-5&1.05E-4&1.02E-4&2.57E-4 \\
0.6& 6.64E-4&8.35E-6&7.53E-5&1.00E-4&9.47E-5&2.44E-4 \\
\noalign{\smallskip}\hline
\end{tabular}
\end{table}

\begin{table}[h]
\caption{The absolute error of Example 5.1 for $\alpha=0.7$ ($0\leq \xi,\eta\leq 0.6$)}
\begin{tabular}{lllllll}
\hline\noalign{\smallskip}
$\xi/\eta$ & $0.1$ & $0.2$ & $0.3$ & $0.4$ & $0.5$ & $0.6$ \\
\noalign{\smallskip}\hline\noalign{\smallskip}
0.1& 2.66E-4&3.77E-5&1.28E-6&1.21E-6&1.63E-5&2.90E-5 \\
0.2& 4.75E-4&6.11E-5&9.13E-6&1.95E-5&4.18E-6&8.61E-5 \\
0.3& 6.26E-4&7.33E-5&2.33E-5&4.07E-5&1.48E-5&1.40E-4 \\
0.4& 7.18E-4&7.78E-5&3.53E-5&5.63E-5&2.95E-5&1.77E-4 \\
0.5& 7.52E-4&7.27E-5&4.66E-5&6.75E-5&4.03E-5&1.97E-4 \\
0.6& 7.25E-4&6.32E-5&5.03E-5&6.69E-5&3.91E-5&1.90E-4 \\
\noalign{\smallskip}\hline
\end{tabular}
\end{table}

\begin{table}[h]
	\caption{The CPU time for Example 5.1 with different $\alpha$ values}
	\begin{tabular}{llll}
		\hline\noalign{\smallskip}
		$\alpha$ & $0.7$ & $0.8$ & $0.9$ \\
		\noalign{\smallskip}\hline\noalign{\smallskip}
		CPU time& 142.40 &148.96 & 150.20 \\
		\noalign{\smallskip}\hline
	\end{tabular}
\end{table}

\noindent$\mathbf{Example ~5.2}$ We take into consideration the following time-fractional Burgers equation (Esen and Tasbozan 2015):
\begin{equation}
D_{\eta}^{\alpha }y-y_{\xi\xi}-yy_{\xi}=f(\xi,\eta)
\end{equation}
\begin{eqnarray*}
0\leq \xi\leq 1, 0\leq \eta\leq 1, \,\ \frac{1}{2}<\alpha\leq 1,
\end{eqnarray*}
\begin{eqnarray*}
f(\xi,\eta)=\frac{4^{\alpha}\eta^{\alpha}\sin(\pi\xi)\Gamma(\alpha+\frac{1}{2})}{\sqrt{\pi}}+\sin(\xi\pi)\pi^2\eta^{2\alpha}-\sin(\xi\pi)\eta^{4\alpha}\cos(\xi\pi)\pi
\end{eqnarray*}
initial and boundary conditions of problem as follow:
\begin{eqnarray}
\left\{
  \begin{array}{ll}
    y(\xi,0)=0, & \\
   y(0,\eta)=y(1,\eta)=0, &
  \end{array}
\right.
\end{eqnarray}
The exact solution of problem is given as:
\begin{equation}
y(\xi,\eta) = \sin(\pi\xi)\eta^{2\alpha}.
\end{equation}
Taking $\xi_i = \frac{i }{p},\,i = 1,\,2,...,p$, $\eta_i = \frac{i }{q},\,i = 1,\,2,...,q$ and $n=p\times q$. In order to demonstrate of the influence and applicability of the method, the absolute error of Example 5.2 is computed and given in Table 5, Table 6 and Table 7 with different values of $\alpha$ for $n=100 (p=q=10)$. The graphics  of approximate solution $y_n \left( \xi,\eta \right)$ are given in Figure 2 with different values of $\alpha=0.7$, $\alpha=0.8$, and $\alpha=0.9$ respectively.

\begin{figure}[ht]
\centerline{\includegraphics[width=2.5000in,height=1.75in]{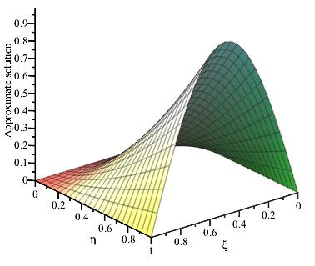}
\includegraphics[width=2.5000in,height=1.75in]{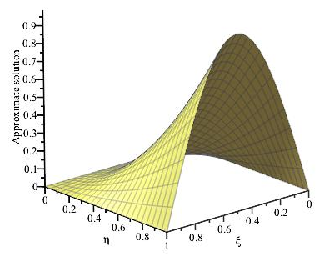}
\includegraphics[width=2.5000in,height=1.75in]{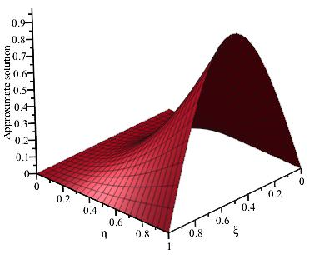}}
\caption{The surface shows approximate solution of Example 5.2 for $n=100$ with $\alpha=0.7$, $\alpha=0.8$ and $\alpha=0.9$ respectively on $0\leq \xi,\eta\leq 1$ }
\label{fig1}
\end{figure}

\begin{table}[h]
\caption{The absolute error of Example 5.2 for $\alpha=0.9$ ($0\leq \xi,\eta\leq 0.6$)}
\begin{tabular}{lllllll}
\hline\noalign{\smallskip}
$\xi/\eta$ & $0.1$ & $0.2$ & $0.3$ & $0.4$ & $0.5$ & $0.6$  \\
\noalign{\smallskip}\hline\noalign{\smallskip}
0.1& 2.50E-3&1.57E-3&5.40E-4&6.93E-5&2.01E-5&9.54E-5 \\
0.2& 4.63E-3&2.70E-3&4.81E-4&7.72E-4&1.37E-3&1.66E-3 \\
0.3& 6.28E-3&3.58E-3&4.21E-4&1.44E-3&2.45E-3&3.05E-3 \\
0.4& 7.32E-3&4.13E-3&3.64E-4&1.89E-3&3.15E-3&3.94E-3 \\
0.5& 7.67E-3&4.30E-3&3.11E-4&2.08E-3&3.43E-3&4.29E-3 \\
0.6& 7.30E-3&4.07E-3&2.64E-4&2.02E-3&3.30E-3&4.10E-3 \\
\noalign{\smallskip}\hline
\end{tabular}
\end{table}

\begin{table}[h]
\caption{The absolute error of Example 5.2 for $\alpha=0.8$ ($0\leq \xi,\eta\leq 0.6$)}
\begin{tabular}{lllllll}
\hline\noalign{\smallskip}
$\xi/\eta$ & $0.1$ & $0.2$ & $0.3$ & $0.4$ & $0.5$ & $0.6$ \\
\noalign{\smallskip}\hline\noalign{\smallskip}
0.1& 2.09E-3&7.98E-4&1.20E-4&8.39E-5&1.77E-4&2.77E-4 \\
0.2& 3.81E-3&1.15E-3&4.53E-4&9.29E-4&1.20E-3&1.52E-3 \\
0.3& 5.16E-3&1.43E-3&9.12E-4&1.74E-3&2.31E-3&2.95E-3 \\
0.4& 6.01E-3&1.59E-3&1.22E-3&2.28E-3&3.04E-3&3.89E-3 \\
0.5& 6.29E-3&1.63E-3&1.36E-3&2.51E-3&3.35E-3&4.28E-3 \\
0.6& 5.98E-3&1.53E-3&1.33E-3&2.45E-3&3.25E-3&4.14E-3 \\
\noalign{\smallskip}\hline
\end{tabular}
\end{table}

\begin{table}[h]
\caption{The absolute error of Example 5.2 for $\alpha=0.7$ ($0\leq \xi,\eta\leq 0.6$)}
\begin{tabular}{lllllll}
\hline\noalign{\smallskip}
$\xi/\eta$ & $0.1$ & $0.2$ & $0.3$ & $0.4$ & $0.5$ & $0.6$ \\
\noalign{\smallskip}\hline\noalign{\smallskip}
0.1& 1.52E-3&2.87E-4&1.56E-4&3.31E-4&4.56E-4&5.71E-4 \\
0.2& 2.69E-3&7.72E-5&5.72E-4&6.69E-4&8.89E-4&1.73E-3 \\
0.3& 3.62E-3&8.15E-5&1.15E-3&1.46E-3&1.96E-3&2.55E-3 \\
0.4& 4.20E-3&1.89E-4&1.53E-3&2.00E-3&2.68E-3&3.46E-3 \\
0.5& 4.39E-3&2.50E-4&1.71E-3&2.26E-3&3.00E-3&3.87E-3 \\
0.6& 4.17E-3&2.67E-4&1.68E-3&2.22E-3&2.94E-3&3.77E-3 \\
\noalign{\smallskip}\hline
\end{tabular}
\end{table}

\begin{table}[h]
	\caption{The CPU time for Example 5.2 with different $\alpha$ values}
	\begin{tabular}{llll}
		\hline\noalign{\smallskip}
		$\alpha$ & $0.7$ & $0.8$ & $0.9$ \\
		\noalign{\smallskip}\hline\noalign{\smallskip}
		CPU time& 536.50 &551.28 & 530.57 \\
		\noalign{\smallskip}\hline
	\end{tabular}
\end{table}

\noindent All computation are made by computer Intel i7-7700HQ CPU-2.8GHz, 16GB RAM on Maple 2017 software.

\section{Conclusion}
In this article, iterative approach of reproducing kernel method has been achievely implemented to find an approximate solution of the time-fractional Burgers equations. The boundedness of linear operator and convergence of iterative process is demonstrated. Numerical results displayed that the present iterative approach is powerful for solving time-fractional Burgers equations with variable and constant coefficient.

\end{document}